\documentclass{article}

\usepackage{amssymb, latexsym}

\usepackage{graphicx}

\usepackage{amsmath}

\def\en t{{{\rm Z}\mkern-5.5mu{\rm Z}}}

\newtheorem{theorem}{Theorem}[section]

\newtheorem{proposition}[theorem]{Proposition}

\newtheorem{remark}[theorem]{Remark}

\begin{document}

\title{\Large\bf  On the $L^p$-estimates of Riesz transforms on forms over complete Riemanian manifolds
}

\author{Xiang-Dong Li\thanks{Research supported by NSFC No. 10971032, Key Laboratory RCSDS, CAS, No. 2008DP173182, and a
Hundred Talents Project of AMSS,
CAS.}\\
\\
{\it\small Academy of Mathematics and Systems Science, Chinese
Academy of Sciences}\\
{\it\small 55, Zhongguancun East Road, Beijing, 100190, P. R. China}\\
{\it\small E-mail: xdli@amt.ac.cn}
}

\maketitle

\begin{center}
\begin{minipage}{120mm}
\begin{center}{\bf Abstract}\end{center} In our previous paper \cite{Li2010}, we
proved a martingale transform representation formula for the Riesz
transforms on forms over complete
Riemannian manifolds, and proved some explicit $L^p$-norm estimates for the Riesz transforms on complete Riemannian manifolds with suitable curvature conditions. In this paper we correct a gap contained in \cite{Li2010} and prove that the main result obtained in \cite{Li2010} on the $L^p$-norm estimates for the Riesz transforms on forms remains valid. Moreover, we prove a time reversal martingale transform representation formula for the Riesz transforms on forms. Finally, we extend our approach and result to the Riesz transforms acting on Euclidean vector bundles over complete Riemannian manifolds with suitable curvature conditions.

\end{minipage}
\end{center}

\section{Introduction}

In our previous paper \cite{Li2010} (Theorem 5.3 p. 507), we obtained the following martingale transform representation formulas for the Riesz transforms on forms over complete Riemannian manifolds:

\begin{eqnarray*}
R_a^1(\square_\phi)\omega(x)=-2\lim\limits_{y\rightarrow +\infty}E_y\left[
\left.  \int_0^\tau e^{a(s-\tau)}M_{\tau, k+1} M_{s, k+1}^{-1}dQ_a\omega(X_s,
B_s)dB_s\right|X_\tau=x\right],\label{R3}
\end{eqnarray*}
\begin{eqnarray*}
R_a^2(\square_\phi)\omega(x)=-2\lim\limits_{y\rightarrow +\infty}E_y\left[
\left.\int_0^\tau e^{a(s-\tau)}M_{\tau, k-1}  M_{s, k-1}^{-1}d_\phi^*Q_a\omega(X_s,
B_s)dB_s\right|X_\tau=x\right].\label{R4}
\end{eqnarray*}
where $R_a^1(\square_{\phi, k})= d(a+\square_{\phi, k})^{-1/2}$ and $R_a^2(\square_{\phi, k})= d_\phi^*(a+\square_{\phi, k})^{-1/2}$.
Recently, Ba\~nuelos and Baudoin \cite{BB} pointed out that,
since $e^{-a\tau}M_{\tau, k\pm 1}$ are not adapted with respect to the filtration  $\mathcal{F}_t=\sigma(X_s: s\in [0, t])$ for $t<\tau$, the above representation formulas should be corrected as follows
\begin{eqnarray}
R_a^1(\square_\phi)\omega(x)=-2\lim\limits_{y\rightarrow +\infty}E_y\left[
\left.  e^{-a\tau}M_{\tau, k+1}  \int_0^\tau e^{as}M_{s, k+1}^{-1}dQ_a\omega(X_s,
B_s)dB_s\right|X_\tau=x\right],\label{R3a}
\end{eqnarray}
\begin{eqnarray}
R_a^2(\square_\phi)\omega(x)=-2\lim\limits_{y\rightarrow +\infty}E_y\left[
\left.e^{-a\tau}M_{\tau, k-1}  \int_0^\tau e^{as}  M_{s, k-1}^{-1}d_\phi^*Q_a\omega(X_s,
B_s)dB_s\right|X_\tau=x\right].\label{R4a}
\end{eqnarray}
Indeed, a careful check of the original proof of Theorem
5.3 in \cite{Li2010} indicates that the correct probabilistic
representation formula of the Riesz transforms $d(a+\square_\phi)^{-1/2}$ and $d_\phi^*(a+\square_\phi)^{-1/2}$ should be  given by
$(\ref{R3a})$ and $(\ref{R4a})$. See Section $2$ below. By the above observation, Ba\~nuelos and
Baudoin \cite{BB} pointed out that
 there is a gap in the proof of the $L^p$-norm
estimates of the Riesz transforms $d(a+\square_\phi)^{-1/2}$ and $d_\phi^*(a+\square_\phi)^{-1/2}$ in
\cite{Li2010} and they proved a new martingale inequality which can be used to correct this gap. In this paper, we correct the above gap and prove that our main result obtained in \cite{Li2010} on the $L^p$-norm estimates of the Riesz transforms on forms  remains valid.
Moreover, we prove a time reversal martingale transform representation formula for the Riesz transforms on forms. Finally, we extend our approach and result to the Riesz transforms acting on Euclidean vector bundles over complete Riemannian manifolds with suitable curvature conditions.

\section{Martingale transform representation formulas}

Let $(M, g)$ be a complete Riemannian manifold, $\nabla$ the
gradient operator on $M$, $\Delta={\rm Tr}\nabla^2$ the covariant Laplace-Beltrami operator on
$M$. Let $\phi\in C^2(M)$, $L=\Delta-\nabla\phi\cdot\nabla$, and
$d\mu=e^{-\phi}dv$, where $dv$ is the standard Riemannian volume
measure on $M$.

Let $d$ be the exterior differential operator, $d^*_\phi$ be its $L^2$-adjoint with respect to the weighted volume measure  $d\mu=e^{-\phi}dv$.
Let $W_k$ be the Weitzenb\"ock curvature operator acting on $k$-forms, and $d\Lambda^k\nabla^2\phi$ be the $k$-linear endomorphism induced by $\nabla^2\phi$ on $\Lambda^kT^*M$.
Let $\square_\phi=dd_\phi^*+d_\phi^*d$ be the Witten Laplacian acting on forms over $(M, g)$ with respect to the weighted volume measure $d\mu=e^{-\phi}dv$.
Recall that the Bochner-Weitzenb\"ock formula reads as
$$\square_{\phi, k}=-(\Delta-\nabla_{\nabla \phi})+W_k+d\Lambda^k\nabla^2\phi.$$
For all $\omega\in C_0^\infty(\Lambda^kT^*M)$, the Poisson integral  $Q_a \omega(x, y)$, also denoted by $\omega_a(x, y)$, is
defined by
$$Q_a\omega(x, y)=e^{-y\sqrt{a+\square_\phi}}\omega(x), \ \ \ \ \forall x\in M, y\geq 0.$$
By \cite{Li2010}, the Riesz transforms associated with the Witten Laplacian  are defined as follows
\begin{eqnarray*}
R_a^1(\square_{\phi, k})&=& d(a+\square_{\phi, k})^{-1/2},\\
R_a^2(\square_{\phi, k})&=& d_\phi^*(a+\square_{\phi, k})^{-1/2}.
\end{eqnarray*}

Let $B_t$ be one dimensional Brownian motion on $\mathbb{R}$
starting from $B_0=y>0$ and with infinitesimal generator ${1\over
2}{d^2\over dy^2}$. Let
$$
\tau=\inf\limits\{t>0: B_t=0\}.$$

Let $X_t$ be the $L$-diffusion process on $M$.
Let $W_t$ be the standard Brownian motion on $\mathbb{R}^n$ such that
$$
dX_t=U_t\circ dW_t-\nabla\phi(X_t)dt,$$
where $U_{t}\in {\rm End}(T_{X_0}M, T_{X_t}M)$ denotes the stochastic parallel transport
along $(X_t)$.
Let $M_{k, t}\in {\rm End}(\Lambda^kT^*_{X_0}M, \Lambda^kT^*_{X_t}M)$ be the solution to the following covariant SDE
along the trajectory of $(X_t)$:
\begin{eqnarray*}
{\nabla M_{t, k}\over \partial t}=-(W_k+d\Lambda^k \nabla^2\phi)(X_t)M_{t, k}, \ \ \ \ M_{0, k}={\rm Id}_{\Lambda^kT^*_{X_0}M}.
\end{eqnarray*}
In the particular case where $W_k+d\Lambda^k \nabla^2\phi=-a$, where $a\geq 0$ is a constant, we have
$$
M_{t, k}=e^{at}U_t, \ \ \ \ \forall t\geq 0.$$

The following results is the correct reformulation of Proposition 5.1 in \cite{Li2010}.

\begin{proposition}\label{pro} For all $\omega\in C_0^\infty(\Lambda^kT^*M)$ and $a\geq 0$, we have
\begin{eqnarray}
\omega(X_\tau)=e^{a\tau}M_{\tau, k}^{*, -1}\omega_a(Z_0)+e^{a\tau}M_{\tau, k}^{*, -1}\int_0^\tau e^{-as}M_{s, k}^*\left(\nabla, {\partial\over \partial y}\right)\omega_a(Z_s)\cdot(U_sdW_s, dB_s).
\label{w1}
\end{eqnarray}
\end{proposition}
{\it Proof}. By It\^o's calculus, we have (see p.504 in \cite{Li2010})
\begin{eqnarray*}
e^{-at}M_{t, k}^*\omega(X_t)=e^{-as}M_{s, k}^*\omega_a(Z_s)+\int_s^t e^{-ar}M_{r, k}^*\left(\nabla, {\partial\over \partial y}\right)\omega_a(Z_r)\cdot(U_rdW_r, dB_r)
\end{eqnarray*}
Taking $s=0$ and $t=\tau$, we obtain Proposition \ref{pro}. \hfill $\square$

The following results is the correct reformulation of  Theorem 5.2 in \cite{Li2010}.

\begin{theorem}\label{theo1} Suppose that $W_k+d\Lambda^k\nabla^2\phi\geq -a$, where $a$ is a non-negative constant. Then, for all  $\omega\in C_0^\infty(\Lambda^kT^*M)$, we have
\begin{eqnarray}
{1\over 2}\omega(x)=\lim\limits_{y\rightarrow \infty}E_y\left[e^{-a\tau}M_{\tau, k}\left.
\int_0^\tau e^{as}M_{s, k}^{-1}{\partial \over \partial y}\omega_a(X_s, B_s)dB_s\right|X_\tau=x\right].\label{w2}
\end{eqnarray}
\end{theorem}
{\it Proof}. The proof is indeed a small modification of the original proof given in \cite{Li2010}. For the completeness of the paper, we give the details here. Let $\eta\in C_0^\infty(\Lambda^kT^*M)$. By $(\ref{w1})$ we have
\begin{eqnarray*}
\eta(X_\tau)=e^{a\tau}M_\tau^{*, -1}\eta_a(Z_0)+e^{a\tau}M_{\tau, k}^{*, -1}\int_0^\tau e^{-as}M_s^*\left(\nabla, {\partial\over \partial y}\right)\eta_a(Z_r)\cdot(U_sdW_s, dB_s).
\end{eqnarray*}
Hence
\begin{eqnarray*}
& &\int_M\left\langle E_y\left[e^{-a\tau}M_{\tau, k}\left.\int_0^\tau e^{as}M_{s, k}^{-1}{\partial \over \partial y}\omega_a(X_s, B_s)dB_s\right|X_\tau=x\right], \eta(x)\right\rangle d\mu(x)\\
& &\ \ \ = E_y\left[e^{-a\tau}M_{\tau, k}\left.\int_0^\tau e^{as}M_{s, k}^{-1}{\partial \over \partial y}\omega_a(X_s, B_s)dB_s, \eta(X_\tau)\right\rangle\right]\\
& &\ \ \ =I_1+I_2,
\end{eqnarray*}
where
\begin{eqnarray*}
I_1&=&E_y\left[\left\langle e^{-a\tau}M_{\tau, k}\int_0^\tau e^{as}M_{s, k}^{-1}{\partial \over \partial y}\omega_a(X_s, B_s)dB_s, e^{a\tau}M_{\tau, k}^{*, -1}\eta_a(X_0, B_0)\right\rangle\right],\\
I_2&=&E_y\left[\left\langle e^{-a\tau}M_{\tau, k}\int_0^\tau e^{as}M_{s, k}^{-1}{\partial \over \partial y}\omega_a(X_s, B_s)dB_s,\right.\right. \\
& & \ \ \ \ \ \ \ \ \ \ \ \ \ \ \   \left.\left. e^{a\tau}M_{\tau, k}^{*, -1}\int_0^\tau e^{-as}M_{s, k}^{*}(\nabla, \partial_y)\eta_a(X_s, B_s)\cdot (U_s dW_s, dB_s)\right\rangle\right].
\end{eqnarray*}
Using the martingale property of the It\^o integral, we have
\begin{eqnarray*}
I_1&=&E_y\left[\left\langle \int_0^\tau e^{as}M_{s, k}^{-1}{\partial \over \partial y}\omega_a(X_s, B_s)dB_s, \eta_a(X_0, B_0)\right\rangle\right]\\
&=&E_y\left[\left\langle E\left[\left.\int_0^\tau e^{as}M_{s, k}^{-1}{\partial \over \partial y}\omega_a(X_s, B_s)dB_s\right|(X_0, B_0)\right], \eta_a(X_0, B_0)\right\rangle\right]\\
&=&0.
\end{eqnarray*}
On the other hand, using the $L^2$-isometry of the It\^o integral, we have
\begin{eqnarray*}
I_2
&=&E_y\left[\int_0^\tau \left\langle  e^{as}M_{s, k}^{-1}{\partial \over \partial y}\omega_a(X_s, B_s), e^{-as}M_{s, k}^{*}{\partial\over \partial y}\eta_a(X_s, B_s)\right\rangle ds\right]\\
&=&E_y\left[\int_0^\tau \left\langle  {\partial \over \partial y}\omega_a(X_s, B_s), {\partial\over \partial y}\eta_a(X_s, B_s)\right\rangle ds\right].
\end{eqnarray*}
The Green function of the background radiation process is given by $2(y\wedge z)$. Hence
\begin{eqnarray*}
& &E_y\left[\int_0^\tau \left\langle  {\partial \over \partial y}\omega_a(X_s, B_s), {\partial\over \partial y}\eta_a(X_s, B_s)\right\rangle ds\right]\\
& &\ \ \ \ \ \ \ \ \ =2\int_M\int_0^\infty (y\wedge z)\left\langle {\partial\over \partial z}\omega_a(x, z),  {\partial\over \partial z}\eta_a(x, z)\right\rangle dzd\mu(x).
\end{eqnarray*}
By spectral decomposition, we have the Littelwood-Paley identity
\begin{eqnarray*}
\lim\limits_{y\rightarrow \infty}\int_M\int_0^\infty (y\wedge z)\left\langle {\partial\over \partial z}\omega_a(x, z),  {\partial\over \partial z}\eta_a(x, z)\right\rangle dzd\mu(x)=\int_M \langle \omega(x), \eta(x)\rangle d\mu(x).
\end{eqnarray*}
Thus
\begin{eqnarray*}
\langle \omega, \eta\rangle_{L^2(\mu)}=2\lim\limits_{y\rightarrow \infty} \int_M \left\langle E_y\left[ e^{-a\tau}M_{\tau, k}\left.
\int_0^\tau e^{as}M_{s, k}^{-1}{\partial\over \partial y}\omega_a(X_s, B_s)dB_s\right|X_\tau=x\right], \eta(x)\right\rangle d\mu(x).
\end{eqnarray*}
This completes the proof of Theorem \ref{theo1}. \hfill $\square$

\medskip

The following martingale transform representation formulas of the Riesz transforms on $k$-forms on complete Riemannian manifolds are the correct reformulations of the ones that we obtained in Theorem $5.3$ in \cite{Li2010}. In the case $k=0$, see \cite{Li2008, Li2013a}.

\begin{theorem}\label{RTF} Under the above notation, for all $\omega\in C_0^\infty(\Lambda^kT^*M)$, we have
\begin{eqnarray}
R_a^1(\square_\phi)\omega(x)=-2\lim\limits_{y\rightarrow +\infty}E_y\left[
\left. e^{-a\tau}M_{\tau, k+1} \int_0^\tau e^{as}M_{s, k+1}^{-1}dQ_a\omega(X_s,
B_s)dB_s\right|X_\tau=x\right],\label{R3}
\end{eqnarray}
\begin{eqnarray}
R_a^2(\square_\phi)\omega(x)=-2\lim\limits_{y\rightarrow +\infty}E_y\left[
\left. e^{-a\tau}M_{\tau, k-1} \int_0^\tau e^{as}M_{s, k-1}^{-1}d_\phi^* Q_a\omega(X_s,
B_s)dB_s\right|X_\tau=x\right].\label{R4}
\end{eqnarray}In particular, in the case where $W_{k+1}+d\Lambda^{k+1}\nabla^2\phi=-a$, we have
\begin{eqnarray}
R_a^1(\square_\phi)\omega(x)=-2\lim\limits_{y\rightarrow
+\infty}E_y\left[\left. U_\tau \int_0^\tau U_s^{-1}dQ_a\omega(X_s,
B_s)dB_s\right|X_\tau=x\right],\label{R5}
\end{eqnarray}
and in the case where $W_{k-1}+d\Lambda^{k-1} \nabla^2\phi=-a$, we have
\begin{eqnarray}
R_a^2(\square_\phi\omega(x)=-2\lim\limits_{y\rightarrow
+\infty}E_y\left[\left. U_\tau \int_0^\tau U_s^{-1}d_\phi^* Q_a\omega(X_s,
B_s)dB_s\right|X_\tau=x\right].\label{R6}
\end{eqnarray}
\end{theorem}
{\it Proof}. The proof is as the same as the one of Theorem $5.3$ in \cite{Li2010}. Indeed,
applying Theorem \ref{theo1} to $R_a^1(\square_\phi)\omega=d(a+\square_{\phi, k})^{-1/2}\omega$,  we have
\begin{eqnarray*}
-{1\over 2}R_a^1(\square_\phi)
\omega(x)&=&\lim\limits_{y\rightarrow \infty}E_y\left[e^{-a\tau}M_{\tau, k+1}
\left.\int_0^\tau e^{as}M_{s, k+1}^{-1}\sqrt{a+\square_{\phi, k+1}}\right.\right.\\
& &\ \ \ \ \ \left.\left. Q_{k+1, a}d(a+\square_{\phi, k})^{-1/2}\omega_a(X_s, B_s)dB_s\right|X_\tau=x\right].
\end{eqnarray*}
Using the commutation formula
$$
d\sqrt{a+\square_{\phi, k}}\omega=\sqrt{a+\square_{\phi, k+1}}d\omega,$$
we obtain
\begin{eqnarray*}
-{1\over 2}R_a^1(\square_\phi)\omega(x)=\lim\limits_{y\rightarrow \infty}E_y\left[e^{-a\tau}M_{\tau, k+1}\left.\int_0^\tau e^{as}
M_{s, k+1}^{-1}dQ_{k, a}\omega_a(X_s, B_s)dB_s\right|X_\tau=x\right].
\end{eqnarray*}
This proves $(\ref{R3})$. Similarly, we can prove  $(\ref{R4})$. Note that, if $W_{k\pm 1}+\nabla^2\phi=-a$, we have $M_{t, k\pm 1}=e^{at}U_{t}$ for all $t\geq 0$. Thus, $(\ref{R5})$ (resp. $(\ref{R6})$) follows from $(\ref{R3})$ (resp. $(\ref{R4})$). \hfill $\square$

\medskip

\begin{remark}{\rm Similarly we have the following martingale transform representation for the Riesz potential on forms.
\begin{eqnarray*}
{1\over 2} (a+\square_\phi)^{-1/2}\omega(x)=-\lim\limits_{y\rightarrow \infty} \left[e^{-a\tau}M_{\tau, k}\left.\int_0^\tau e^{as}
M_{s, k}^{-1}\omega_a(X_s, B_s)dB_s\right|X_\tau=x\right].
\end{eqnarray*}
In particular, under the condition $W_k+d\Lambda^k\nabla^2\phi\geq 0$, we have
\begin{eqnarray*}
{1\over 2} \square_\phi^{-1/2}\omega(x)=-\lim\limits_{y\rightarrow \infty} \left[M_{\tau, k}\left.\int_0^\tau
M_{s, k}^{-1}\omega(X_s, B_s)dB_s\right|X_\tau=x\right],
\end{eqnarray*}
where $\omega(x, y)=e^{-y\sqrt{\square_\phi}}\omega(x)$ denotes the Poisson semigroup generated by $\square_\phi$ on $L^2(\Lambda^kT^*M, \mu)$.}
\end{remark}
\section{The $L^p$-norm estimate}

In this section we correct a gap contained in \cite{Li2010} and prove that our main result obtained in \cite{Li2010} on the $L^p$-norm estimates of the Riesz transforms on forms remains valid. When $p=2$, we have the following

\begin{proposition}\label{RTF2} For all $a\geq 0$ and $\omega\in C_0^\infty(\Lambda^kT^*M)$, we have
\begin{eqnarray*}
\|d(a+\square_\phi)^{-1/2}\omega\|_2&\leq& \|\omega\|_2,\\
\|d_\phi^*(a+\square_\phi)^{-1/2}\omega\|_2&\leq&\|\omega\|_2.
\end{eqnarray*}
\end{proposition}
{\it Proof}. By Gaffney's integration by part, for all $\omega\in C_0^\infty(\Lambda^kT^*M)$, we have
\begin{eqnarray*}
\langle\langle (a+\square_\phi)\omega, \omega\rangle\rangle=a\|\omega\|_2^2+\|d\omega\|_2^2+\|d_\phi^*\omega\|_2^2.
\end{eqnarray*}
Since $a+\square_\phi$ is non-negative symmetric operator on $L^2(\Lambda^kT^*M, \mu)$, we get
\begin{eqnarray*}
\|d\omega\|_2^2+\|d_\phi^*\omega\|_2^2+a\|\omega\|_2^2=\|\sqrt{a+\square_\phi}\omega\|_2^2.
\end{eqnarray*}
This implies that
\begin{eqnarray*}
\|d(a+\square_\phi)^{-1/2}\omega\|_2^2+\|d_\phi^*(a+\square_\phi)^{-1/2}\omega\|_2^2\leq \|\omega\|_2^2.
\end{eqnarray*}
The proof of Proposition \ref{RTF2} is completed. \hfill $\square$

\medskip

The following result is the restatement of the main result (i.e., Theorem 1.6) in \cite{Li2010}.

\begin{theorem}\label{EST2} Let $M$ be a complete Riemannian manifold, and $\phi\in
C^2(M)$. Suppose that there exists a constant $a\geq 0$ such that
$$W_k+d\Lambda^k \nabla^2\phi\geq -a, \ \ \ {\rm and}\ \ \
W_{k+1}+d\Lambda^{k+1} \nabla^2\phi\geq -a.$$ Then, there exists a
constant $C_k>0$ depending only on $k$ such that for all $p>1$,
\begin{eqnarray} \|d(a+\square_\phi)^{-1/2}\omega\|_p\leq C_k(p^*-1)^{-3/2}\|\omega\|_p, \ \ \
\forall \omega\in C_0^\infty(\Lambda^kT^*M). \ \label{f1}
\end{eqnarray}
In particular, if $W_k+d\Lambda^k\nabla^2\phi\geq 0$ and $W_{k+1}+d\Lambda^{k+1}\nabla^2\phi\geq 0$, then the Riesz transform $d\square_{\phi}^{-1/2}$ is bounded in $L^p$ for all $p>1$, and there exists a
constant $C_k>0$ depending only on $k$ such that for all $p>1$,
\begin{eqnarray} \|d\square_\phi^{-1/2}\omega\|_p\leq C_k(p^*-1)^{-3/2}\|\omega\|_p, \ \ \
\forall \omega\in C_0^\infty(\Lambda^kT^*M). \label{f2}
\end{eqnarray}
\end{theorem}
{\it Proof}. By Theorem \ref{RTF}, Fatou's lemma, and using the $L^p$-contractivity of the conditional expectation, for any $1<p<\infty$, we have
\begin{eqnarray*}
& &\|d(a+\square_\phi)^{-1/2}\omega\|_p^p\\
&=&2^p\int_M \lim\limits_{y\rightarrow \infty}\left|E_y\left[e^{-a\tau}M_{\tau, k+1}\left.
\int_0^\tau e^{as} M_{s, k+1}^{-1}dQ_{a, k}\omega(X_s, B_s)dB_s\right|X_\tau=x\right]\right|^pd\mu(x)\\
&\leq&2^p \lim\inf\limits_{y\rightarrow \infty}\int_M E_y\left[\left|e^{-a\tau}M_{\tau, k+1}\left.
\int_0^\tau e^{as} M_{s, k+1}^{-1}dQ_{a, k}\omega(X_s, B_s)dB_s\right|^p\right|X_\tau=x\right]d\mu(x)\\
&=&2^p\lim\inf\limits_{y\rightarrow \infty}E_y\left[ \left| e^{-a\tau}M_{\tau, k+1}
\int_0^\tau e^{as} M_{s, k+1}^{-1}dQ_{a, k}\omega(X_s, B_s)dB_s\right|^p\right].
\end{eqnarray*}
Recall that, see p. 509-p. 510 in \cite{Li2010},  there is an $(n+1)\times (n+1)$ operator valued matrix such that
$$
d\omega(x, y)=A\overline{\nabla}\omega(x, y),$$
where $\overline{\nabla}=(\nabla, \partial_y)$. Moreover, $\|A\|_{\rm op}$ is a finite number depending only on $k$. In view of this, we have
\begin{eqnarray}
\|d(a+\square_\phi)^{-1/2}\omega\|_p\leq 2\lim\inf\limits_{y\rightarrow \infty}\left\| e^{-a\tau}M_{\tau, k+1}
\int_0^\tau e^{as} M_{s, k+1}^{-1}A\overline{\nabla}Q_{a, k}\omega(X_s, B_s)\cdot (U_sdW_s, dB_s)\right\|_p.\label{TTT}
\end{eqnarray}
Let
\begin{eqnarray*}
I_y= e^{-a \tau}M_{\tau, k+1}\int_0^\tau e^{as}M_{s,
k+1}^{-1}A\overline{\nabla}Q_{a, k}\omega(X_s, B_s)\cdot (U_sdW_s,
dB_s),
\end{eqnarray*}
and
\begin{eqnarray*}
J_y=\left\{\int_0^\tau |\overline{\nabla}Q_{a, k}\omega(X_s,
B_s)|^2ds\right\}^{1/2}.
\end{eqnarray*}

By Theorem 2.6 due to Ba\~nuelos and Baudoin in \cite{BB}, under the assumption $W_k+d\Lambda^k \nabla^2\phi
\geq -a$, we can prove that
\begin{eqnarray}
\|I_y\|_p\leq 3\sqrt{p(2p-1)}\|A\|_{\rm op}\|J_y\|_p.\label{IJ}
\end{eqnarray}
Moreover, by Proposition $6.2$ in our previous paper \cite{Li2010}, we have
\begin{eqnarray*}
\|J_y\|_p\leq B_p\|\omega\|_p,
\end{eqnarray*}
where $B_p=(2p)^{1/2}(p-1)^{-3/2}$ for $p\in (1, 2)$, $B_p=1$ for $p=2$, and
$B_p={p\over \sqrt{2(p-2)}}$ if $p>2$. Combining this with
$(\ref{IJ})$, for all $1<p<2$, we can prove that
\begin{eqnarray*}
\|R_a^1(\square_\phi)\omega\|_p&\leq &2\lim\inf\limits_{y\rightarrow \infty}\|I_y\|_p\\
&\leq& 6\sqrt{2}\|A\|_{\rm op}p(2p-1)^{1/2}(p-1)^{-3/2}\|\omega\|_p\\
&\leq& 12\sqrt{6}\|A\|_{\rm op}(p-1)^{-3/2}\|\omega\|_p,
\end{eqnarray*}
and for $p>2$, we have
\begin{eqnarray*}
\|R_a^1(\square_\phi)\omega\|_p&\leq& 2\lim\inf\limits_{y\rightarrow \infty}\|I_y\|_p\\
&\leq& 3\sqrt{2}\|A\|_{\rm op}p^{3/2}(2p-1)^{1/2}(p-2)^{-1/2}\|\omega\|_p\\
&\leq& 6\|A\|_{\rm op}(p-1)^{3/2}(1+O(1/p))\|\omega\|_p.
\end{eqnarray*}
This implies the desired $L^p$-norm estimate for the Riesz transform $d(a+\square_\phi)^{-1/2}$. \hfill $\square$

\begin{remark}{\rm The above proof corrects a gap in the proof of Theorem 1.6 in \cite{Li2010} (p. 510 line 4 to line 5 in \cite{Li2010}),  where we used the
Burkholder-Davies-Gundy inequality to derive that
\begin{eqnarray*}
\|I_y\|_p\leq C_p\left\| \left\{ \int_0^\tau |e^{a(s-\tau)}M_{\tau,
k+1}M_{s, k+1}^{-1}A|^2 |\overline{\nabla}Q_{a, k}\omega(X_s,
B_s)|^2ds \right\}^{1/2} \right\|_p,\label{Iy}
\end{eqnarray*}
where $C_p$ is a constant depending only on $p$. However,  as $e^{-a\tau}M_{\tau, k\pm 1}$ are not adapted with respect to the filtration  $\mathcal{F}_t=\sigma(X_s: s\in [0, t])$ for $t<\tau$, one cannot use the Burkholder-Davis-Gundy
inequality in the above way, except that $e^{-a\tau}M_{\tau, k\pm 1}$ is independent of $(X_s, s\in [0, \tau])$, which only happens in the case where  $W_{k+1}+d\Lambda^{k+1}\nabla^2\phi\equiv -a$.}
\end{remark}

\medskip

\begin{theorem}\label{EST3} Let $M$ be a complete Riemannian manifold, and $\phi\in
C^2(M)$. Suppose that there exists a constant $a\geq 0$ such that
$$W_{k}+d\Lambda^k \nabla^2\phi\geq -a, \ \ \ {\rm and}\ \ \
W_{k-1}+d\Lambda^{k-1} \nabla^2\phi\geq -a.$$ Then, there exists a
constant $C_k>0$ depending only on $k$ such that for all $p>1$,
\begin{eqnarray} \|d_\phi^*(a+\square_\phi)^{-1/2}\omega\|_p\leq C_k(p^*-1)^{3/2}\|\omega\|_p, \ \ \
\forall \omega\in C_0^\infty(\Lambda^kT^*M). \ \label{f1}
\end{eqnarray}
In particular, if $W_k+d\Lambda^k\nabla^2\phi\geq 0$ and $W_{k-1}+d\Lambda^{k-1}\nabla^2\phi\geq 0$, then the Riesz transform $d\square_{\phi}^{-1/2}$ is bounded in $L^p$ for all $p>1$. More precisely, there exists a
constant $C_k>0$ depending only on $k$ such that for all $p>1$,
\begin{eqnarray} \|d_\phi^*\square_\phi^{-1/2}\omega\|_p\leq C_k(p^*-1)^{3/2}\|\omega\|_p, \ \ \
\forall \omega\in C_0^\infty(\Lambda^kT^*M). \label{f2}
\end{eqnarray}
\end{theorem}
{\it Proof}. By duality argument as used in \cite{Li2010}, we can derive Theorem \ref{EST3} from Theorem \ref{EST2}. \hfill $\square$

\section{Case of constant curvature}

In the particular case where $W_{k+1}+d\Lambda^{k+1}\nabla^2\phi\equiv -a$,
we have
\begin{eqnarray*}
d(a+\square_\phi)^{-1/2}\omega(x)=-2\lim\limits_{y\rightarrow
+\infty}E_y\left[\left. U_\tau \int_0^\tau U_s^{-1}A\overline{\nabla}Q_{a, k}\omega(X_s,
B_s)dB_s\right|X_\tau=x\right].
\end{eqnarray*}
By the same argument as used in the proof of $(\ref{TTT})$,
\begin{eqnarray*}
\|d(a+\square_\phi)^{-1/2}\omega\|_p\leq 2\lim\inf\limits_{y\rightarrow \infty}\left\| U_\tau \int_0^\tau U_s^{-1}A\overline{\nabla}Q_{a, k}\omega(X_s,
B_s)\cdot (U_SdW_s, dB_s)\right\|_p.
\end{eqnarray*}
Hence
\begin{eqnarray*}
\|d(a+\square_\phi)^{-1/2}\omega\|_p\leq 2\lim\inf\limits_{y\rightarrow \infty}\left\| \int_0^\tau U_s^{-1}A\overline{\nabla}Q_{a, k}\omega(X_s,
B_s)\cdot(U_sdW_s, dB_s)\right\|_p.\label{TT}
\end{eqnarray*}
By Burkholder's sharp $L^p$-inequality for subordination of martingale transforms \cite{Bk} we have
\begin{eqnarray*}
\|d(a+\square_\phi)^{-1/2}\omega\|_p\leq 2\|A\|_{\rm op}(p^*-1)\lim\inf\limits_{y\rightarrow \infty}\left\| \int_0^\tau U_s^{-1}\overline{\nabla}Q_{a, k}\omega(X_s,
B_s)\cdot(U_sdW_s, dB_s)\right\|_p.
\end{eqnarray*}As was pointed out in Remark 6.5 in \cite{Li2010}, only if $W_k+d\Lambda^k \nabla^2\phi=-a$, we can obtain
\begin{eqnarray*}
\left\| \int_0^\tau U_s^{-1}\overline{\nabla}Q_{a, k}\omega(X_s,
B_s)\cdot(U_sdW_s, dB_s)\right\|_p\leq \|\omega\|_p.
\end{eqnarray*}
That is to say, only if $W_k+d\Lambda^k\nabla^2\phi=-a$ and $W_{k+1}+d\Lambda^{k+1}\nabla^2\phi=-a$, which happens in the case where $M$ is a flat Riemannian manifold and $\nabla^2\phi\equiv 0$, hence $a=0$,
we can obtain
$\|d\square_\phi^{-1/2}\omega\|_p\leq 2\|A\|_{\rm op}(p^*-1)\|\omega\|_p$.

\begin{remark}{\rm  In view of Theorem \ref{EST2}, for all $p>1$, the upper bound $C_k(p^*-1)^{3/2}$ appeared in Theorem 1.6 in \cite{Li2010} remains valid, but the upper bound $C_k(p^*-1)$ appeared in Theorem 1.7 and Theorem 1.8 in \cite{Li2010} should be replaced by $C_k(p^*-1)^{3/2}$.}
\end{remark}

\section{Time reversal martingale transformation representation formula for the Riesz transfroms}

In this section, we prove a time reversal martingale transformation representation formula
for the Riesz transforms on forms on complete Riemannian
manifolds.

\medskip
First, we have the following time reversal martingale transformation representation formula for forms.

\begin{theorem}\label{theo1-b} Let $\widehat{X}_t=X_{\tau-t}$, and $\widehat{B}_t=B_{\tau-t}$, $t\in [0,
\tau]$. Let $\widehat{M}_{t, k}$ be the solution to the covariant SDE
\begin{eqnarray*}
{\nabla\over \partial t}\widehat
{M}_{t, k}&=&-\widehat{M}_t(W_k+d\Lambda^k \nabla^2\phi)(\widehat{X}_t),\\
\widehat{M}_{0, k}&=&{\rm Id}_{\Lambda^k T^*_{\widehat X_0}M}.
\end{eqnarray*}
For any $\omega\in C_0^\infty(\Lambda^k T^*M)$, let $\omega_a(x,
y)=e^{-y\sqrt{a+\square_\phi}}\omega(x)$, $\forall x\in M, y\geq 0$.
Then, for a.s. $x\in M$,
\begin{eqnarray*}
{1\over 2}\omega(x)=\lim\limits_{y\rightarrow +\infty}E_y\left[\left.\widehat Z_\tau\right|\widehat{X}_0=x\right],
\end{eqnarray*}
where
\begin{eqnarray*}
\widehat Z_\tau=\int_0^\tau
e^{-at}\widehat{M_{t, k}}\partial_y\omega_a(\widehat{X}_t,
\widehat{B}_t)d\widehat{B}_t-\int_0^\tau e^{-at}
\widehat{M}_{t, k}\partial_y^2 \omega(\widehat{X}_t, \widehat{B}_t)dt.
\end{eqnarray*}
\end{theorem}
{\it Proof}. The proof is similarly to the one of Theorem 5.1 in \cite{Li2013a}.  \hfill $\square$

By Theorem \ref{theo1-b}, we can prove the following time reversal martingale
transformation representation formula for the Riesz transforms on
complete Riemannian manifolds.

\begin{theorem}  For any $\omega\in C_0^\infty(M, \Lambda^kT^*M)$, we have
\begin{eqnarray*}
R_a^1(\square_{\phi, k})\omega(x)&=&-2\lim\limits_{y\rightarrow +\infty}E_y\left[
\left. \widehat Z_{\tau, k+1} \right|\widehat{X}_0=x\right],\\
R_a^2(\square_{\phi, k})\omega(x)&=&-2\lim\limits_{y\rightarrow +\infty}E_y\left[
\left. \widetilde Z_{\tau, k-1} \right|\widehat{X}_0=x\right],
\end{eqnarray*}
where
\begin{eqnarray*}
\widehat Z_{\tau, k+1}&=&\int_0^\tau e^{-as}\widehat M_{s, k+1} dQ_{k, a}\omega(\widehat{X}_s,
\widehat{B}_s)d\widehat{B}_s-\int_0^\tau e^{-as}\widehat M_{s, k+1}\partial_ydQ_{k, a}\omega(\widehat{X}_s, \widehat{B}_s)ds,\\
\widetilde Z_{\tau, k-1}&=&\int_0^\tau e^{-as}\widehat M_{s, k-1} d_\phi^* Q_{k, a}\omega(\widehat{X}_s,
\widehat{B}_s)d\widehat{B}_s-\int_0^\tau e^{-as}\widehat M_{s, k-1}\partial_yd_\phi^*Q_{k, a}\omega(\widehat{X}_s, \widehat{B}_s)ds.
\end{eqnarray*}
\end{theorem}

\begin{remark} {\rm
As noticed in \cite{Li2008}, there exists a standard one dimensional Brownian motion $\beta_t$ such that
\begin{eqnarray*}
d\widehat{B}_t=d\beta_t+{dt\over \widehat{B}_t}, \ \ \ t\in (0, \tau].
\end{eqnarray*}}
\end{remark}

\section{Riesz transforms on Euclidean vector bundles}

In this section we extend our approach and result to the Riesz transforms acting on Euclidean vector bundles over complete Riemannian manifolds.

Let $M$ be a complete Riemannian manifold, $E$ a Riemannian vector bundle over $M$. Let $\nabla^E$ be a metric preserving connection on $E$.
Let $F=\Lambda^\cdot T^*M\otimes E$, and define
$$\nabla^F=\nabla^{\Lambda^\cdot T^*M}\otimes 1_E+1_{\Lambda^\cdot T^*M}\otimes \nabla^E.$$
The De Rham operator acting on $C^\infty(M, F)$ is defined by
$$
d^F=\sum\limits_{i=1}^n e_i^*\wedge \nabla_{e_i}^F,
$$
where $(e_1, \ldots, e_n)$ is a orthonormal basis at any point $x\in M$, and $(e_1^*, \ldots, e_n^*)$ is its dual.

The curvature of $\nabla^E$ is defined by
\begin{eqnarray*}
R^E=({\nabla^E})^2.
\end{eqnarray*}
Suppose that $E$ is an Euclidean vector bundle with flat connection, i.e, $R^E=0$. Then
\begin{eqnarray*}
(d^F)^2=0.
\end{eqnarray*}

Let $\phi\in C^2(M)$, $\mu=e^{-\phi}dv$. Let $d_\phi^{F*}$ be the $L^2$-adjoint of $d^F$ with respect to $\mu$. We have
\begin{eqnarray*}
(d_\phi^{F*})^2=0.
\end{eqnarray*}
The Witten Laplacian acting on $C_0^\infty(M, F)$ is defined by
\begin{eqnarray*}
\square_{F, \phi}=d^Fd^{F*}_\phi+d^{F*}_\phi d^F.
\end{eqnarray*}
The heat semigroup and the Poisson semigroup generated by $\square_{F, \phi}$ are denoted by $P_t\omega(x)=e^{-t\square_{F, \phi}}\omega(x)$ and  $Q_a\omega(x, y)=e^{-y\sqrt{a+\square_{F, \phi}}}\omega(x)$ respectively. The Bochner-Weitzenb\"ock formula holds
\begin{eqnarray*}
\square_{F, \phi}=-\Delta_{F, \phi}+W_{F, \phi},
\end{eqnarray*}
where $\Delta_{F, \phi}={\rm Tr}(\nabla^F)^2-\nabla_{\nabla\phi}^F$, and $W_{F, \phi}=W+d\Lambda^\cdot \nabla^2\phi$.

Let $X_t$ be the $L$-diffusion process on $M$. Let $M_{k, t}\in {\rm End}(\Lambda^kT^*_{X_0}M\otimes E, \Lambda^kT^*_{X_t}M\otimes E)$ be the solution to the following covariant SDE
along the trajectory of $(X_t)$:
\begin{eqnarray*}
{\nabla M_{t, k}\over \partial t}=-(W_k+d\Lambda^k \nabla^2\phi)(X_t)M_{t, k}, \ \ \ \ M_{0, k}={\rm Id}_{\Lambda^kT^*_{X_0}M\otimes E}.
\end{eqnarray*}

We have the following results on the quantitative $L^p$-estimates of the Riesz transforms on Euclidean vector bundles over complete Riemannian manifolds.

\begin{theorem}\label{Main-A} Let $M$ be a complete Riemannian manifold, $E$ be an Euclidean vector bundle over $M$, and $\phi\in C^2(M)$. Then, for all $\omega\in C_0^\infty(M, \Lambda^kT^*M\otimes E)$ and for all $\mu$-a.s. $x\in M$, we have
\begin{eqnarray*}
d^F(a+\square_{F, \phi})^{-1/2}\omega(x)=-2\lim\limits_{y\rightarrow \infty}E_y\left[
\left.  e^{-a\tau}M_{\tau, k+1}  \int_0^\tau e^{as}M_{s, k+1}^{-1}d^F
Q_a\omega(X_s, B_s)dB_s\right|X_\tau=x\right].
\end{eqnarray*}
Suppose that $W_{i}+d\Lambda^{i}\nabla^2\phi\geq -a$, $i=k,\ k+1$. Then, for all $p>1$ and for all $\omega\in L^p(\Lambda^{k}T^*M\otimes E, \mu)$, we have
\begin{eqnarray*}
\|d^F(a+\square_{F, \phi})^{-1/2}\omega\|\leq C_p\|A\|\|\omega\|_p
\end{eqnarray*}
where $\|A\|$ is the operator norm of $A\in {\rm End}(F, F)$ is such that $d^F\omega=A\nabla\omega$ and depends only on $k$, $C_p$ is a constant depending only on $p$, more precisely, $C_p=\|A\|^{-1}$ for $p=2$, $C_p=O(p^*-1)^{3/2}$ for $p\rightarrow 1$ and $p\rightarrow \infty$.
\end{theorem}
{\it Proof}. The proof is as the same as the one of Theorem \ref{RTF} and Theorem \ref{EST2}. \hfill $\square$

\begin{theorem}\label{Main-B} Let $M$ be a complete Riemannian manifold, $E$ be an Euclidean vector bundle over $M$, and $\phi\in C^2(M)$. Then, for all $\omega\in C_0^\infty(M, \Lambda^kT^*M\otimes E)$ and for all $\mu$-a.s. $x\in M$, we have
\begin{eqnarray*}
d^{F*}_{\phi}(a+\square_{F, \phi})^{-1/2}\omega(x)=-2\lim\limits_{y\rightarrow \infty}E_y\left[
\left.  e^{-a\tau}M_{\tau, k-1}  \int_0^\tau e^{as}M_{s, k-1}^{-1}d^{F*}_{\phi}
Q_a\omega(X_s, B_s)dB_s\right|X_\tau=x\right],
\end{eqnarray*}
Suppose that $W_{i}+d\Lambda^{i}\nabla^2\phi\geq -a$, $i=k,\ k-1$. Then, for all $p>1$ and for all $\omega\in L^p(\Lambda^{k}T^*M\otimes E, \mu)$, we have
\begin{eqnarray*}
\|d^{F*}_\phi (a+\square_{F, \phi})^{-1/2}\omega\|\leq C_p\|A\|\|\omega\|_p
\end{eqnarray*}
where $C_p$ is a constant depending only on $p$, more precisely, $C_p=\|A\|^{-1}$ for $p=2$, and $C_p=O(p^*-1)^{3/2}$ for $p\rightarrow 1$ and $p\rightarrow \infty$.
\end{theorem}
{\it Proof}. By duality argument, we can derive Theorem \ref{Main-B} from Theorem \ref{Main-A}. \hfill $\square$

\medskip

To end this paper, let us mention that, in a forthcoming paper \cite{Li2013b},  we will prove a martingale transform representation formula for the Riesz transforms associated with the Dirac operator acting on Hermitian vector bundles over complete Riemannian manifolds and for the Riesz transforms associated with the $\bar\partial$-operator acting on holomorphic Hermitian vector bundles over complete K\"ahler manifolds. By the same argument as used in this paper and in \cite{Li2013a}, we can prove some explicit dimension free $L^p$-norm estimates of these Riesz transforms on complete Riemannian or K\"ahler manifolds with suitable curvature conditions. See also \cite{Li2010b}.

\medskip

\noindent{\bf Acknowledgement}.   \ \ I would
like to thank R. Ba\~nuelos and F. Baudoin
for their interests on my previous works and for pointing out the gap contained in \cite{Li2010}
which has been addressed in this paper. I would like also to thank Yong Liu and Songzi Li for helpful discussions and encouragements.

\end{document}